\def\softd{{\leavevmode\setbox1=\hbox{d}%
\hbox to 1.05\wd1{d\kern-0.4ex{\char039}\hss}}}%cstocs
\def\softt{{\leavevmode\setbox1=\hbox{t}%
\hbox to \wd1{t\kern-0.6ex{\char039}\hss}}}%cstocs
\def\softl{l\kern-0.45ex\raise0.1ex\hbox{'}\kern-0.10ex}%cstocs
\def\softL{L\kern-0.8ex\raise0.1ex\hbox{'}\kern0.1ex}%cstocs
\documentstyle[a4,psfig]{article}
\begin{document}
\thispagestyle{empty}
\vspace*{40mm}
\begin{center}
        {\Large {\bf{THE MODELLING AND ANALYSIS OF FRACTIONAL-ORDER
                     CONTROL SYSTEMS IN FREQUENCY DOMAIN}}} \\
         \vspace{10mm}
         \normalsize
         {I. PETR\'A\v{S}$^{\mbox{1}}$, \softL{}.    		DOR\v{C}\'AK$^{\mbox{1}}$,
          P. O'LEARY$^{\mbox{2}}$, B. M. VINAGRE$^{\mbox{3}}$,
          I. PODLUBN\'Y$^{\mbox{1}}$} \\
               \vspace{1mm}
               $^{\mbox{1}}$~{ Department of Informatics and Process Control  \\
                    BERG Faculty, Technical University of Ko\v{s}ice   \\
                    B. N\v{e}mcovej 3, 042 00 Ko\v{s}ice, Slovak Republic  \\
                    phone:      (+42195) 6025172                      \\
                    e-mail: {\it \{petras, dorcak, podlbn\}@tuke.sk}} \\
               \vspace{1mm}
               $^{\mbox{2}}$~{ Department of Automation, University of Leoben   \\
                    Peter - Tunner strasse 27, A-8700 Leoben, Austria  \\
                    phone:      (+43) 3842/4029031                      \\
                    e-mail: {\it oleary@unileoben.ac.at}} \\
               \vspace{1mm}
               $^{\mbox{3}}$~{ Department of Electronic and Electromechanic Engineering  \\
                    Industrial Engineering School, University of Extremadura   \\
                    Avda. De Elvas, s/n., 06071-Badajoz, Spain  \\
                    phone:      (+34) 924289600                      \\
                    e-mail: {\it bvinagre@unex.es}} \\

\end{center}
\vskip 5mm
\begin{abstract}
\hspace{-6.5mm}
This paper deals with fractional-order controlled systems and
fractional-order controllers in the frequency domain. The mathematical
description by fractional transfer functions and properties of
these systems are presented. The new ways for modelling of
fractional-order systems are illustrated with a~numerical example
and obtained results are discussed in conclusion.
\end{abstract}

\hspace{2.5mm}

\section*{1. INTRODUCTION}

Fractional calculus is a generalization of integration and derivation
to non-integer order fundamental operator $_{a}D^{\alpha}_{t}$,
where $a$ and $t$ are the limits of the operation. The two definitions
used for the general fractional differintegral
are the Gr\"unwald definition and the Riemann-Liouville definition
\cite{Oldham}.

The idea of fractional calculus has been known since the development
of the regular calculus, with the first reference
probably being associated with Leibniz and L'Hospital in 1695.
Fractional calculus was used for modelling  of physical systems,
but we can find only few works dealing with the application of
this mathematical tool in control theory (e.g. \cite{Axtell,
Petras1, Oustaloup, Podlubny2}).

The aim of this paper is to show, how by using the fractional
calculus,
we can obtain a more general structure for the classical $PID$
controller, for controlled systems with memory and hereditary
behaviour, and how to model fractional-order systems
in the frequency domain.
On the other hand, we can analyse the fractional-order systems
and specify the conditions of stability in Bode's and
Nyquist's frequency response.

\newpage
\thispagestyle{empty}
\section*{2. FRACTIONAL-ORDER CONTROL CIRCUIT}

\begin{figure}[h]
\centerline{\psfig{file=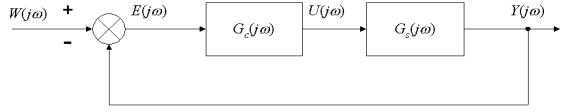,width=8.7cm}}
\caption{Feed - back control loop}
\label{F1}
\end{figure}

We will be studying the control system shown in Fig.\ref{F1},
where $G_{c}(j\omega)$ is the controller transfer function
and $G_{s}(j\omega)$ is the controlled system transfer function.

\subsection*{2.1 Fractional-order controlled system}

The fractional controlled system will be represented with fractional
model with frequency transfer function given by the following expression:
\begin{equation}\label{Gs}
      G_{s}(j\omega) =\frac{Y(j\omega)}{U(j\omega)}=
                \frac{b_{m}(j\omega)^{\alpha_{m}} + \ldots
     + b_{1}(j\omega)^{\alpha_{1}} + b_{0}(j\omega)^{\alpha_{0}}}
     {a_{n}(j\omega)^{\beta_{n}} + \ldots
     + a_{1}(j\omega)^{\beta_{1}} + a_{0}(j\omega)^{\beta_{0}}}
     =\frac{\sum_{k=0}^{m}b_k(j\omega)^{\alpha_k}}
        {\sum_{k=0}^{n}a_k(j\omega)^{\beta_k}},
\end{equation}
where  $\beta_k, \alpha_k$  $(k = 0, 1, 2, \dots, n)$ are generally real
numbers,
$\beta_{n} > \ldots > \beta_{1} > \beta_{0}$,
$\alpha_{m} > \ldots > \alpha_{1} > \alpha_{0}$ and
$a_k, b_k$ $(k= 0, 1, \dots, n)$ are arbitrary constants.

Identification  methods \cite{Ludovic, Oustaloup, Podlubny3,
Vinagre} for determination of the coefficients
$a_k, b_k$ and  $\alpha_k, \beta_k$ $(k= 0, 1, \dots, n)$
were developed, based on minimisation of the difference between
the calculated and experimentally measured values. The optimisation
consists in minimising the square norm of the difference between
the measured frequency response $F(\omega)$ and the frequency
response of the model $G_{s}(j\omega)$. The quadratic criterion
is given by the expression:
\begin{equation}
Q=\sum_{m=0}^{M} W^2(\omega_m)
|F(\omega_m)-G_{s}(j\omega_m)|^2,
\end{equation}
where $W(\omega_m)$ is a weighting function and $M$ is the number
of measured values of frequencies.

\subsection*{2.2 Fractional order controller}

The fractional $PI^{\lambda}D^{\delta}$ controller will be represented
by frequency transfer function given in the following expression:
\begin{equation}\label{Gr}
      G_{c}(j\omega)=\frac{U(j\omega)}{E(j\omega)}=
      K + \frac{T_i}{(j\omega)^{\lambda}}+T_d(j\omega)^{\delta},
\end{equation}
where $\lambda$ and $\delta$ are arbitrary real numbers
$(\lambda, \delta \geq 0)$, $K$ is the proportional constant,
$T_i$ is the integration constant and $T_d$ is the derivation
constant.

Taking $\lambda=1$ and $\delta=1$, we obtain a~classic $PID$
controller. If $\lambda = 0$ and/or $T_i=0$, we obtain
a~$PD^{\delta}$ controller, etc.
All these types of controllers are particular cases of the
$PI^{\lambda}D^{\delta}$ controller, which is more flexible and
gives an opportunity to better adjust the dynamical properties
of the fractional-order control system.

The $PI^{\lambda}D^{\delta}$ controller with complex zeros and
poles located anywhere in the left-hand $s$-plane
may be rewritten as
\begin{equation}
  G_{c}(j\omega)= C
\frac{\left((j\omega)/\omega_n\right)^{\delta+\lambda}
+(2\xi(j\omega)^{\lambda})/\omega_n+1}
{(j\omega)^{\lambda}},
\end{equation}
where $C$ is a gain, $\xi$  is the dimensionless damping ratio and
$\omega_n$ is the natural frequency. Normally, we choose $0.9 >
\xi >0.7$. When $\xi=1$, the condition is called critical
damping \cite{Dorf}.

The tuning of $PI^{\lambda}D^{\delta}$ controller parameters is
determined according to the given requirements. These requirements
are, for example, the stability measure, the accuracy of the
regulation process, dynamical properties etc.
One of the methods being developed is the method (modification
of roots locus method) of dominant roots \cite{Petras1}, based on
the given stability measure and the damping measure of the control
circuit. Another possible way to obtain the controller parameters
is using the tuning formula, based on gain and phase margins
specifications.

\section*{3. FREQUENCY CHARACTERISTICS }

The graphical interpretation of the frequency transfer function
for different values of the angular
frequency  $\omega$  in  the range   $\omega \in \langle 0,
\infty \rangle$  is called the frequency characteristic.

The practical meaning of the frequency characteristics is in
the determination of the stability of control circuits from the
behavior of the frequency characteristic. Stability is an asymptotic
qualitative criterion of the quality of the control circuit and is
the primary and necessary condition for correct functioning of
every control circuit. It is difficult to evaluate the stability of
a fractional-order control circuit by examining its characteristic
equation either by finding its dominant roots or by algebraic methods.
The methods given are suitable for integer-order systems. In our case
we use Bode's and Nyquist's frequency characteristic
to investigate the stability \cite{Petras2}.

For the frequency transfer of the open control circuit
$G_o(j\omega)$ with respect to (\ref{Gs}) and (\ref{Gr}), we have
\begin{equation}\label{Go}
    G_o(j\omega)=\frac{Y(j\omega)}{W(j\omega)} =
\left ( \frac{T_i + K (j\omega)^{\lambda} +
T_d(j\omega)^{\delta+\lambda}}{(j\omega)^{\lambda}} \right )
        \frac{\sum_{k=0}^{m}b_k(j\omega)^{\alpha_k}}
        {\sum_{k=0}^{n}a_k(j\omega)^{\beta_k}}.
\end{equation}

It is known from the investigation theory of stability of regulation
circuits that the system is stable if the roots of the characteristic
equation are negative or have negative real parts if they are complex
conjugate. This means they are located left of the imaginary axis of
the complex $s$-plane of the roots. In the case of frequency methods
of evaluating the stability we transform the plane of the roots
$s$ into the complex plane $G_o(j\omega)$ while the transformation
is realised according to the transfer function of the open circuit.
During the transformation all the roots of the characteristic equation
are mapped from the $s$-plane into the critical point $(-1, i0)$
in the plane $G_o(j\omega)$. The mapping of the $s$-plane into plane
$G_o(j\omega)$ is conformal, that is, the directions and locations
of points in the $s$-plane  are preserved also in the plane
$G_o(j\omega)$.

In the Nyquist frequency characteristic to find out the stability
of the control circuit, it is necessary to investigate the
behavior of the curve $G_o(j\omega)$ for
$\omega \in \langle 0,\infty \rangle$
relative to the critical point $(-1; i0)$. Based on the above it
follows that also the fractional-order control circuit is stable,
if the frequency characteristic passes on the right side of the
critical point, when going in the direction of increasing values
$\omega$. If the frequency characteristic passes on the left of
the critical point, the circuit is non-stable. The passing of the
frequency characteristic through the critical point means the circuit
is on the border of stability.

\section*{4. ILLUSTRATIVE EXAMPLE}

We give in this section an example of modelling the
stable dynamical system by using fractional calculus in
frequency domain. The fractional-order control system consists of
the real controlled system with the frequency transfer function
\begin{equation}\label{sys}
        G_s(j\omega)=\frac{1}
                          {0.8 (j\omega)^{2.2} + 0.5 (j\omega)^{0.9} + 1}
\end{equation}
and the fractional-order $PD^{\delta}$ controller, designed on the
stability measure $S_t=2.0$ and damping measure $\xi=0.4$, with
the frequency transfer function in the form:
\begin{equation}\label{cont}
%        G_c(j\omega)=50.15 + 392.82 (j\omega)^{0.9} + 24.98(j\omega)^{0.9}.
         G_c(j\omega)=50.0 + 5.326 (j\omega)^{1.286}.
\end{equation}

%\begin{figure}[h]
    \vskip 0.5mm
     \centerline{\psfig{file=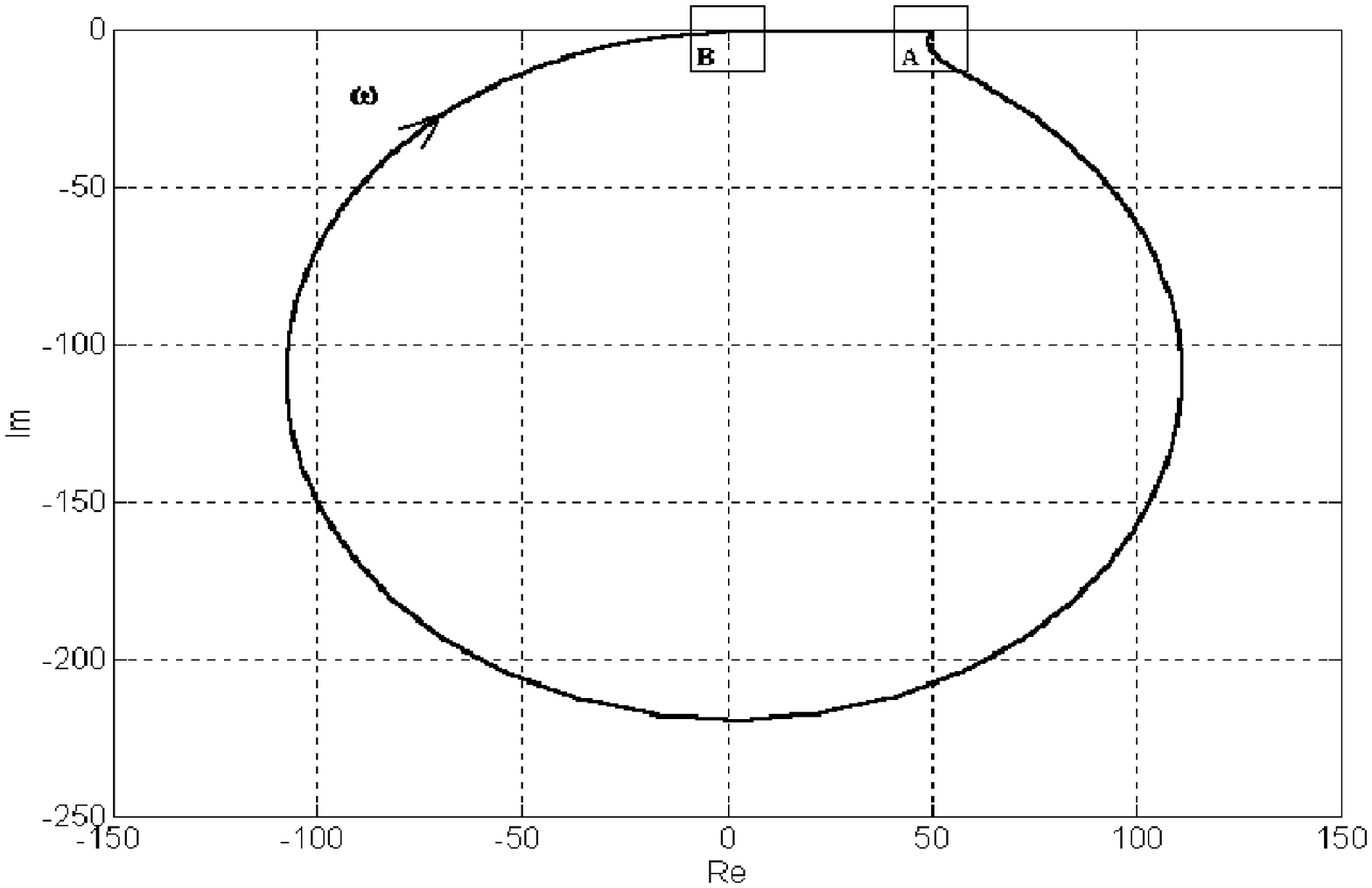,width=7cm}
		 \psfig{file=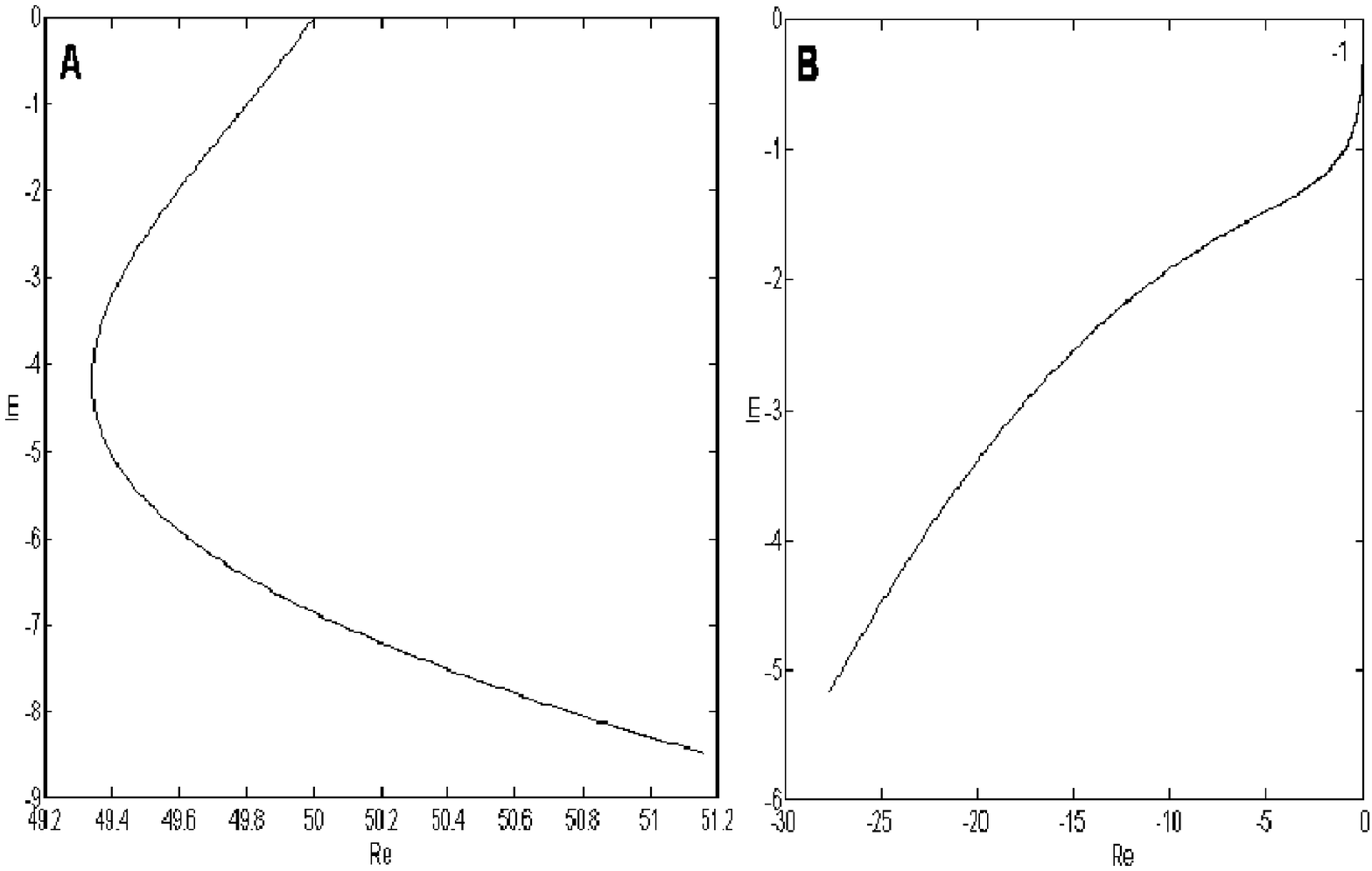,width=7cm}} %100 percent
     \centerline{Figure 2:\,Nyquist plot (and details) for (\ref{sys}) and (\ref{cont})}
%\end{figure}
\thispagestyle{empty}
\newpage
\thispagestyle{empty}
\vspace*{-7mm}
%\begin{figure}[h]
     \centerline{\psfig{file=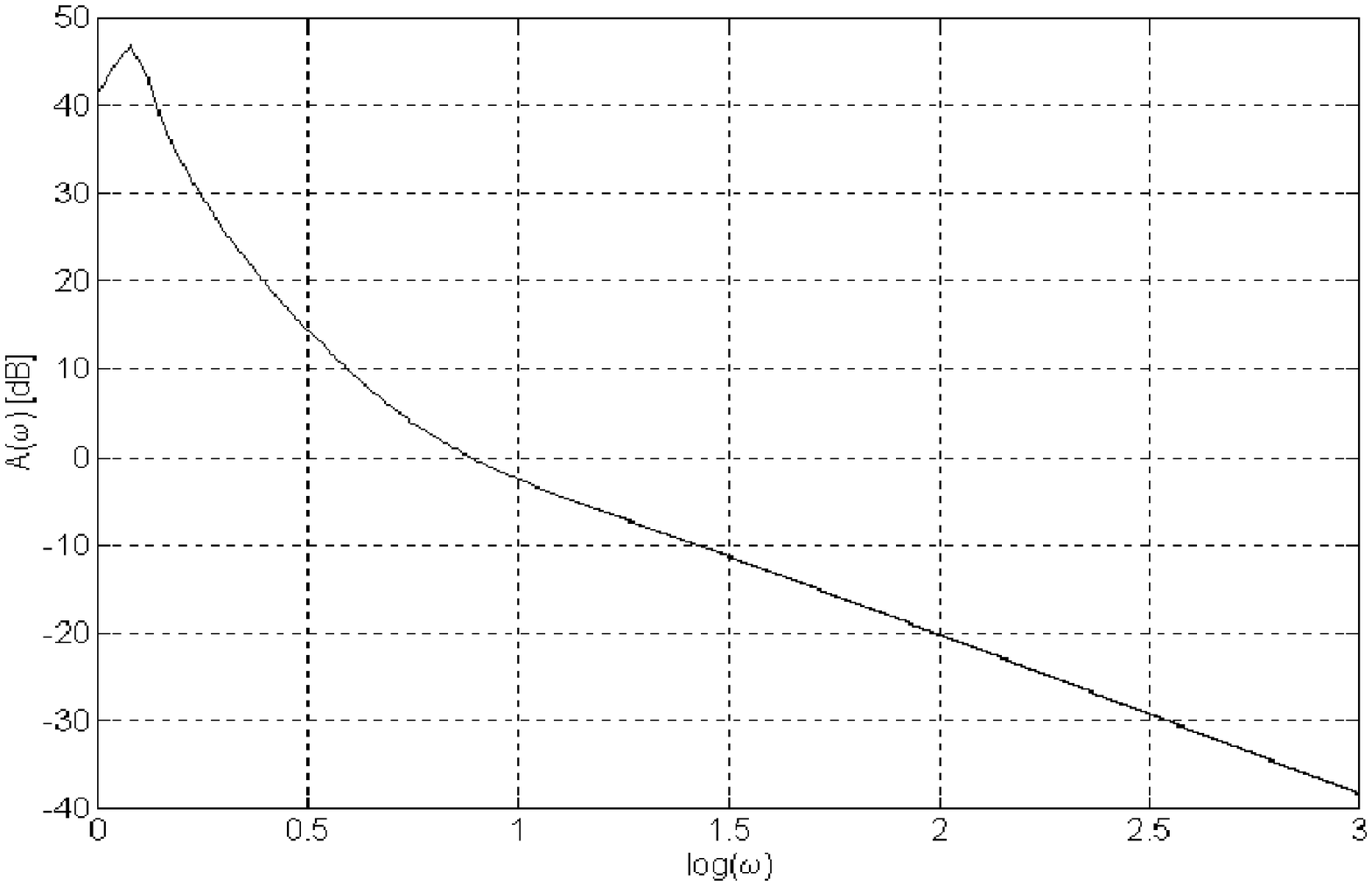,width=7cm}
		 \psfig{file=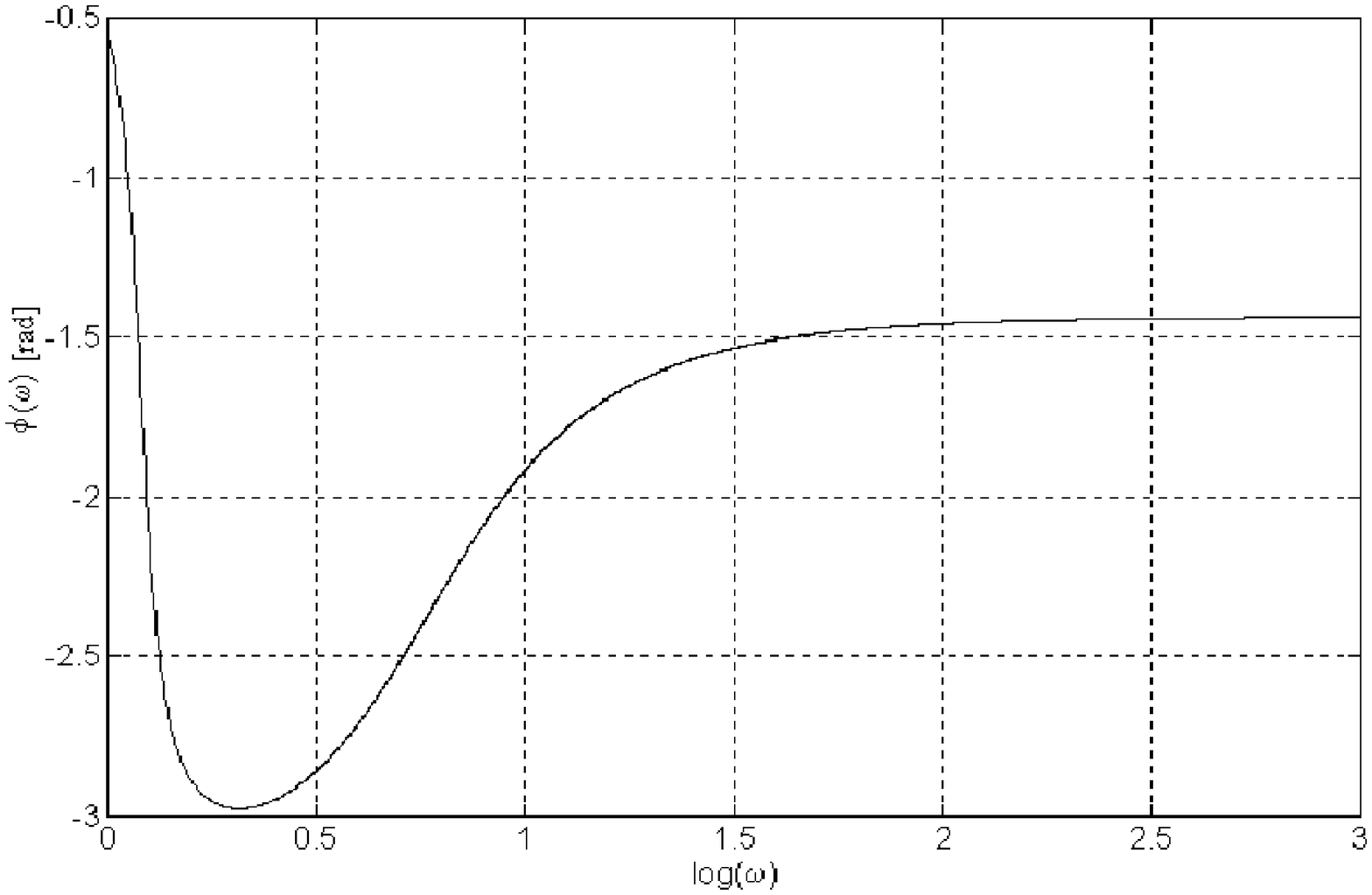,width=7cm}} %100 percent
     \centerline{Figure 3:\,Bode plots for (\ref{sys}) and (\ref{cont})}
     \vskip 5 mm
%\end{figure}

In Fig. 2 and Fig. 3 are showed the frequency characteristics of
the fractional-order control system consists of the fractional-order
controlled object (\ref{sys}) and the fractional-order
controller (\ref{cont}).

\section*{5. CONCLUSION}

The above methods make it possible to model and analyse
fractional-order control systems in the frequency domain.
The stability of a fractional-order control system can be
investigated via the behavior of the frequency plot (see Fig. 2 and
Fig. 3). These properties of the fractional-order systems can be used in the controller parameters design.
The results of this and previous works also show that fractional-order
controllers are robust.
This can even lead to qualitatively different dynamical phenomena
in control circuits.

\section*{ACKNOWLEDGEMENTS}

This work was partially supported by grant {\it \,Fractional-order Dynamical Systems and Controllers:
Discrete and Frequency-Domain Models and Algorithms} from
Austrian Institute for East and South-East Europe
and partially supported by grant VEGA 1/7098/20 from the
Slovak Agency for Science.

\end{document}